\documentclass{article}
\usepackage{amsmath}
\usepackage{amsfonts}
\usepackage{amssymb}
\usepackage{euscript}
\input epsf

\newtheorem{thm}{Theorem}[section]

\newtheorem{example}{Example}[section]

\def\dj{d\kern-0.4em\char"16\kern-0.1em}

\title{DEFORMED MITTAG--LEFFLER POLYNOMIALS
\thanks{Supported by
Ministry of Sci. \& Techn. Rep. Serbia, the projects  No. 144023.}
}
\author{
{\sc Miomir S. Stankovi\'c}\\
Faculty of Ocupational Safety\\
[2mm]
{\sc Sladjana D. Marinkovi\'c}\\
Faculty of Electronic Engineering\\
[2mm]
{\sc Predrag M. Rajkovi\'c}\\
Faculty of Mechanical Engineering\\
[2mm]
{\bf   University of Ni\v s,\ Serbia}
}

\date{}

\begin{document}
\maketitle

\medskip

{\it Abstract.} The starting point of this paper are the Mittag--Leffler polynomials introduced by H. Bateman in \cite{Bateman1}. Based on generalized integer powers of real numbers and deformed exponential function, we introduce deformed Mittag--Leffler polynomials defined by appropriate generating function. We investigate their recurrence relations, differential properties and orthogonality. Since they have all zeros on imaginary axes, we also consider real polynomials with real zeros associated to them.

\medskip





\section{Introduction}

The development in various disciplines of modern physics leads to including some deformed and generalized versions of the exponential functions (see \cite{Tsallis2} and \cite{Kaniadakis1}). In that direction, we defined two variable deformed exponential function in \cite{DefExpFun}. Here, we will use it to define a class of polynomials, related to Mittag--Leffler polynomials.

Let us remind that Mittag--Leffler polynomials $g_n(y)$ are introduced by H. Bateman (see \cite{Bateman1}, \cite{Bateman2}) as coefficients in expansion
$$
(1+x)^y(1-x)^{-y}=\sum_{n=0}^\infty g_n(y)x^n\qquad (|x|<1).
$$
Although they are known for a long time, a few new papers considering them, appear recently (see, for example, \cite{He} and \cite{Luzon}).

They satisfy recurrence relations
$$
g_n(y+1)-g_{n-1}(y+1)=g_n(y)+g_{n-1}(y),
$$
\begin{equation*}
(n+1)g_{n+1}(y)-2yg_n(y)+(n-1)g_{n-1}(y)=0,
\end{equation*}
both  with initial values
\begin{equation*}
 g_0(y)=1,\quad g_1(y)=2y,
\end{equation*}
and orthogonality relation
\begin{equation}\label{MLPOrtho}
\int_{-\infty}^{+\infty}g_n(-iy)g_m(iy)\frac{dy}{y\sinh \pi y}=\frac{2}{n}\delta_{mn}\qquad(n,m\in\mathbb N).
\end{equation}

Notice that the corresponding monic sequence
\begin{equation*}\label{MLPMonic}
\hat g_n(y) = \frac{(n)!}{2^{n}}\ g_n(y)  \qquad (n\in\mathbb N_0)
\end{equation*}
has the exponential generating function
\begin{equation*}\label{GenFMLPMonic}
\Bigl(\frac{2+x}{2-x}\Bigr)^y = \sum_{n=0}^\infty \hat g_n(y)\frac{x^n}{n!}  \ .
\end{equation*}

Let $h\in\mathbb R\setminus\{0\}$. We define generalized integer powers of real numbers  \cite{Riordan}
\begin{eqnarray*}
z^{(0,h)}=z^{[0,h]}=1,\qquad
z^{(n,h)}=\prod_{k=0}^{n-1}(z-kh)\ ,
\quad
z^{[n,h]}=\prod_{k=0}^{n-1}(z+kh)\quad (n\in\mathbb N),
\end{eqnarray*}

Recall the $h$--difference operator  \cite{Jagerman}
\begin{equation}\label{forwOp}
\Delta_{z,h} f(z)=\frac{f(z+h) - f(z)}{h}
=\frac1{h}\bigl(E_h-I\bigr)f(z),
\end{equation}
where $I$ is identical and $E_h$ is shift--operator.
This $h$--difference operator are linear. However, product rules for them are:
\begin{equation}
\Delta_{z,h}\bigl(f(z)g(z)\bigr)=f(z+h)\Delta_{z,h}g(z)+\Delta_{z,h}f(z)g(z).\label{Delta-prod}
\end{equation}
Its acting on integer generalized powers is given by:
$$
\Delta_{z,h}\,z^{(n,h)}=nz^{(n-1,h)},\qquad
\Delta_{z,h}\,z^{[n,h]}=n(z+h)^{[n-1,h]}.
$$

In \cite{DefExpFun}, it was defined deformed exponential function
\begin{equation*}
e_h(x,y)=(1+hx)^{y/h} \qquad (x\in \mathbb C\setminus\{-1/h\},\ y\in \mathbb R).
\end{equation*}
Function $e_h(x,y)$ keeps some of basic properties of exponential function. For $y\in \mathbb R$, the following holds:
$$
\aligned
e_h(x,y)&>0\qquad (x<-1/h\quad \text{for}\ h<0\quad \text{or}\quad x>-1/h\quad \text{for}\ h>0), \\
e_h(0,y)&=e_h(x,0)=1.
\endaligned
$$
If $h$ exchanges the sign, we have
\begin{equation}
e_{-h}(x,y)=e_h(-x,-y)\qquad (x\ne 1/h).\label{e-h}
\end{equation}
The additional property is kept only in regard to second variable:
$$
e_h(x,y_1)e_h(x,y_2)=e_h(x,y_1+y_2).
$$

Deformed exponential functions can be represented as expansions:
\begin{eqnarray}
e_h(x,y)&=&\sum_{n=0}^\infty \dfrac{1}{n!}\ x^n y^{(n,h)}\qquad(|hx|<1),\label{sum_eh} \\
e_{-h}(x,y)&=&\sum_{n=0}^\infty \dfrac{1}{n!}\ x^n y^{[n,h]}\qquad(|hx|<1).\label{sum_e-h} 
\end{eqnarray}
If we apply $h$--difference operators on deformed exponential functions, we get:
\begin{eqnarray}
\Delta_{y,h}\ e_h(x,y)&=&\sum_{n=1}^\infty \dfrac{1}{n!}\ x^n n y^{(n-1,h)}=x\ e_h(x,y),\label{Delta_eh}\\
\Delta_{y,h}\ e_{-h}(x,y)&=&\sum_{n=1}^\infty \dfrac{1}{n!}\ x^n n (y+h)^{[n-1,h]}=x\ e_{-h}(x,y+h).\label{Delta_e-h}
\end{eqnarray}
An interesting differential property of this function (see \cite{DefExpFun}) is given as
\begin{equation}
\Bigl((1+hx)\frac{\partial}{\partial x}\Bigr)\ e_h(x,y)=y\ e_h(x,y).\label{partial_eh}
\end{equation}

\section{Deformed Mittag--Leffler polynomials}

Generating function of Mittag--Leffler polynomials can be recognized as
$$
G(x,y)=(1+x)^y(1-x)^{-y}=e_1(x,y)e_{-1}(x,y).
$$
For $h\in\mathbb R\setminus \{0\}$ we can define {\it deformed Mittag--Leffler polynomials} as coefficients in expansion
\begin{equation}\label{MLpol1}
G_h(x,y)=e_h(x,y)e_{-h}(x,y)=\sum_{n=0}^\infty g_n^{(h)}(y)x^n
\end{equation}
With respect to (\ref{sum_eh}) and (\ref{sum_e-h}), we have
$$
\aligned
e_h(x,y)e_{-h}(x,y)&=\sum_{n=0}^\infty \frac{y^{(n,h)}}{n!}x^n \sum_{m=0}^\infty \frac{y^{[m,h]}}{m!}x^m \\
&=\sum_{n=0}^\infty \sum_{m=0}^n \frac{y^{(m,h)}y^{[n-m,h]}}{m!(n-m)!}x^n.
\endaligned
$$
Hence,
$$
g_n^{(h)}(y)=\sum_{m=0}^n \frac{y^{(m,h)}y^{[n-m,h]}}{m!(n-m)!}
=\frac{1}{n!}\sum_{m=0}^n \binom nm y^{(m,h)}y^{[n-m,h]}\qquad(n\in\mathbb N_0).
$$

According to (\ref{e-h}), we have
$$
G_h(x,y)=G_h(-x,-y)=\sum_{n=0}^\infty g_n^{(h)}(-y)(-1)^n x^n,
$$
and, consequently,
$$
g_n^{(h)}(-y)=(-1)^n g_n^{(h)}(y).
$$
Also, from \ $G_h(x,y)=G_{-h}(x,y)$ we get
$$
g_n^{(-h)}(y)=g_n^{(h)}(y).
$$

Let us derive some recurrence relations for polynomials $g_n^{(h)}(y)$.

\begin{thm}\label{thm1}
The successive members of sequence $\{g_n^{(h)}(y)\}_{n\in\mathbb N_0}$ satisfy the three--term recurrence relation
\begin{eqnarray*}
(n+1)g_{n+1}^{(h)}(y)-2yg_n^{(h)}(y)-h^2(n-1)g_{n-1}^{(h)}(y)&=&0\quad (n\ge 2),\\
\quad g_0^{(h)}(y)=1,\qquad g_1^{(h)}(y)&=&2y.
\nonumber
\end{eqnarray*}
\end{thm}

\noindent{\it Proof.} Firstly, we have
$$
\frac{\partial}{\partial x}G_h(x,y)=
\left(\frac{\partial}{\partial x}e_h(x,y)\right)e_{-h}(x,y)+e_h(x,y)\left(\frac{\partial}{\partial x}e_{-h}(x,y)\right).
$$
Multiplying by $(1-hx)(1+hx)$, with respect to the differential property (\ref{partial_eh}) of $e_h(x,y)$, we obtain:
$$
\aligned
(1-h^2 x^2)&\frac{\partial}{\partial x}G_h(x,y) \\
&=(1-hx)ye_h(x,y)e_{-h}(x,y)+(1+hx)e_h(x,y)ye_{-h}(x,y)\\
&=2yG_h(x,y).
\endaligned
$$
Using (\ref{MLpol1}) and comparing the coefficients in the series we obtain the recurrence relation. $\ \Box$

\begin{example}\rm
The first members of the sequence $\{g_n^{(h)}(y)\}_{n\in\mathbb N_0}$ are:
$$
g_0^{(h)}(y)= 1,\quad g_1^{(h)}(y)=2y,\quad g_2^{(h)}(y)=2y^2,
\quad g_3^{(h)}(y)=\frac23 y(2y^2+h^2),
$$
$$
g_4^{(h)}(y)=\frac23 y^2(y^2+2h^2),\quad
g_5^{(h)}(y)=\frac2{15}y(2y^4+10h^2y^2+3 h^4),
$$
\end{example}

Let us remark the useful relation
\begin{equation}\label{deltaGx}
\frac{\partial}{\partial x}G_h(x,y)=\frac{2y}{1-h^2 x^2}G_h(x,y)
\end{equation}
obtained in the proof of Theorem \ref{thm1}.

\begin{thm} For polynomials $g_n^{(h)}(y)$ the following relation is valid:
$$
g_{n}^{(h)}(y+h)-g_{n}^{(h)}(y)=h\bigl(g_{n-1}^{(h)}(y+h)+g_{n-1}^{(h)}(y)\bigr)\qquad (n\in\mathbb N).
$$
\end{thm}

\noindent{\it Proof.}
According to (\ref{Delta-prod}), (\ref{Delta_eh}) and (\ref{Delta_e-h}) we have
$$
\aligned
\Delta_{y,h}G_h(x,y)&=\Delta_{y,h}\bigl(e_h(x,y)e_{-h}(x,y)\bigr)\\
&=e_h(x,y+h)\Delta_{y,h}e_{-h}(x,y)+e_{-h}(x,y)\Delta_{y,h}e_h(x,y)\\
&=x\bigl(e_h(x,y+h)e_{-h}(x,y+h)+e_h(x,y)e_{-h}(x,y)\bigr)\\
&=x\bigl(G_h(x,y+h)+G_h(x,y)\bigr).
\endaligned
$$
Using (\ref{MLpol1}), because of linearity of operator $\Delta_{y,h}$ the last equation becomes
$$
\aligned
\sum_{n=0}^\infty \Delta_{y,h}g_n^{(h)}(y)x^n&=x\sum_{n=0}^\infty \bigl(g_n^{(h)}(y+h)+g_n^{(h)}(y)\bigr)x^n\\
&=\sum_{n=1}^\infty \bigl(g_{n-1}^{(h)}(y+h)+g_{n-1}^{(h)}(y)\bigr)x^n.
\endaligned
$$
With respect to (\ref{forwOp}), we obtain required relation. $\ \Box$

\begin{thm} The polynomial $g_n^{(h)}(y)$ can be represented over hypergeometric function as
\begin{equation}\label{hyper}
g_n^{(h)}(y)=2yh^{n-1}\;{}_2F_1\Bigl({{1-n,\ 1-y/h}\atop{2}}\Bigm| 2\Bigr)\qquad (n\in\mathbb N).
\end{equation}
\end{thm}

\noindent{\it Proof.} It follows from the previous theorem by mathematical induction. Namely, its equivalent form is summation formula
$$\aligned
g_n^{(h)}(y)
&=2yh^{n-1}\ \sum_{k=0}^{n-1} \frac{(1-n)_k (1-y/h)_k}{(2)_k}\frac{2^k}{k!}\\
&=2yh^{n-1}\ \sum_{k=0}^{n-1} (1-y/h)_k\frac{(-1)^k \ 2^k}{k!(k+1)!} \prod_{j=1}^k (n-j) .
\endaligned$$
It is true for $n=1$. We will suppose that it is valid for all $k\le n$ and check validity of three term recurrence relation. Then
$$
\aligned
(n+1)g_{n+1}^{(h)}(y)&-2yg_n^{(h)}(y)-h^2(n-1)g_{n-1}^{(h)}(y)\\
&=2yh^{n}\Bigl((n+1)\sum_{k=0}^{n} (1-y/h)_k\frac{(-1)^k \ 2^k}{k!(k+1)!} \prod_{j=1}^k (n+1-j)\\
&\qquad\qquad+2\sum_{k=0}^{n-1} (1-y/h)_k(-y/h)\frac{(-1)^k \ 2^k}{k!(k+1)!} \prod_{j=1}^k (n-j) \\
&\qquad\qquad -(n-1)\sum_{k=0}^{n-2} (1-y/h)_k \ \frac{(-1)^k \ 2^k}{k!(k+1)!} \prod_{j=1}^k (n-1-j)\Bigr).
\endaligned
$$
Since
$$
(1-y/h)_{k+1}=(1-y/h)_k(-y/h)+(k+1)(1-y/h)_k,
$$
the coefficient beside $(1-y/h)_k$ $(k=0,1,\ldots, n)$ in the recurrence relation is
$$
\frac{(-1)^k 2^k}{k!(k+1)!} \Bigl( (n+1)n-(n-k)(n-k-1)-k(k+1)-2(k+1)(n-k)\Bigr)\prod_{j=1}^{k-1}(n-j)=0.\ \Box
$$

Let us show that the polynomials of sequence $\{g_n^{(h)}(y)\}_{n\in\mathbb N_0}$ are orthogonal.

\begin{thm} For the members of sequence of polynomials $\{g_n^{(h)}(y)\}_{n\in\mathbb N_0}$ \linebreak
the following orthogonality relation is valid:
\begin{equation}\label{orth}
\int_{-\infty}^{\infty}g_n^{(h)}(iy)g_m^{(h)}(-iy)\frac{dy}{y\sinh(\pi y/h)}=\frac{2h^{2n-2}}{n}\;\delta_{mn}\qquad(m,n\in\mathbb N).
\end{equation}
\end{thm}

\noindent{\it Proof.} Using representation (\ref{hyper}) and the integral representation of the hypergeometric function, we can write
$$
\aligned
g_n^{(h)}(y)&=2yh^{n-1}\frac{\Gamma(2)}{\Gamma(1-y/h)\Gamma(1+y/h)}\int_0^1 u^{-y/h}(1-u)^{y/h}(1-2u)^{n-1}du\\
&=\frac{h^n}{\pi}\sin\frac{\pi y}{h}\int_{-1}^1 u^{n-1}(1+u)^{y/h}(1-u)^{-y/h}du,
\endaligned
$$
or, by changing of variable $u=\tanh (ht/2)$,
$$
g_n^{(h)}(y)=\frac{h^n}{\pi}\sin\frac{\pi y}{h}
\int_{-\infty}^{\infty} \Bigl(\tanh\frac{ht}2\Bigr)^{n}\frac{e^{ty}}{\sinh ht}\;dt.
$$
Hence,
$$
\aligned
g_n^{(h)}(iy)&=\frac{ih^n}{\pi}\sinh\frac{\pi y}{h}
\int_{-\infty}^{\infty} e^{ity}\Bigl(\tanh\frac{ht}2\Bigr)^{n}\frac{dt}{\sinh ht}\\
&=ih^n\sqrt{\frac2{\pi}}\;\sinh\frac{\pi y}{h}\;\mathcal F\left(\Bigl(\tanh\frac{ht}2\Bigr)^{n}\frac{1}{\sinh ht}\right),
\endaligned
$$
where $\varphi(t)\mapsto \Phi(y)=\mathcal F\left(\varphi(t)\right)$ denotes the Fourier transform.
Applying the inverse Fourier transform, we get
$$
\Bigl(\tanh\frac{ht}2\Bigr)^{n}\frac{1}{\sinh ht}
=\frac{1}{ih^n}\sqrt{\frac{\pi}{2}}\ \mathcal F^{-1}\left(\frac{g_n^{(h)}(iy)}{\sinh (\pi y/h)}\right),
$$
i.e.,
\begin{equation}\label{tanh^n}
\Bigl(\tanh\frac{ht}2\Bigr)^{n}=\frac{1}{2ih^n}\sinh{ht}\int_{-\infty}^{\infty}e^{-ity}\frac{g_n^{(h)}(iy)}{\sinh (\pi y/h)}\;dy.
\end{equation}
Further, according to (\ref{deltaGx}), we recognize
$$
\aligned
\sinh{ht}\ e^{-ity}&=\frac{2\tanh(ht/2)}{1-\tanh^2(ht/2)} \left(\frac{1+\tanh(ht/2)}{1-\tanh(ht/2)}\right)^{-iy/h}\\
&=\frac{ih}{y}\tanh\frac{ht}{2}\ \frac{\partial}{\partial x}G_h(x,-iy)\Bigm|_{x=\frac{1}{h}\tanh(ht/2)}.
\endaligned
$$
Because of
$$
\frac{\partial}{\partial x}G_h(x,-iy)\Bigm|_{x=\frac{1}{h}\tanh(ht/2)}
=\sum_{m=1}^\infty \frac{m}{h^{m-1}}g_m^{(h)}(-iy)\left(\tanh\frac{ht}{2}\right)^{m-1},
$$
we obtain
$$
\sinh{ht}\ e^{-ity}=\frac{ih}{y}\sum_{m=1}^\infty \frac{m}{h^{m-1}}g_m^{(h)}(-iy)\left(\tanh\frac{ht}{2})\right)^{m}.
$$
Substituting last equation in (\ref{tanh^n}), we have
$$
\Bigl(\tanh\frac{ht}2\Bigr)^{n}
=\frac{1}{2h^{n-1}}\sum_{m=1}^\infty \frac{m}{h^{m-1}}\left(\tanh\frac{ht}{2}\right)^{m}
\int_{-\infty}^{\infty}\frac{g_n^{(h)}(iy)g_m^{(h)}(-iy)}{y\sinh (\pi y/h)}\;dy.
$$
Comparing the coefficients by $\tanh(ht/2)$, we get required relation.$\ \Box$

The monic sequence $\{\hat g_n^{(h)}(y)\}_{n\in\mathbb N_0}$ is related with $\{g_n^{(h)}(y)\}_{n\in\mathbb N_0}$ by
$$
g_n^{(h)}(y)= \frac{2^n}{n!} \hat g_{n}^{(h)}(y).
$$
They satisfy the three--term recurrence relation
$$
\aligned
&\hat g_{n+1}^{(h)}(y) = y\ \hat g_n^{(h)}(y)+h^2\frac{n(n-1)}4 \ \hat g_{n-1}^{(h)}(y)  \qquad (n\in\mathbb N),\\
&\quad \hat g_0^{(h)}(y)=1,\qquad \hat g_1^{(h)}(y)=y,
\endaligned
$$
and have generating function given by
$$
\hat G_h(x,y)=e_h(x/2,y)e_{-h}(x/2,y)=\sum_{n=0}^\infty \hat g_n^{(h)}(y)\frac{x^n}{n!}\ .
$$

\section{The polynomials associated to Mittag--Leffler and deformed Mittag--Leffler polynomials}

Let us consider {\it modified Mittag-Leffler polynomials} defined by
\begin{equation}\label{ModMLP}
\varphi_n(y) = \frac{g_{n+1}(iy)}{i^{n+1}\ y}\qquad (n\in\mathbb N_0).
\end{equation}

\begin{thm}
The successive members of sequence $\{\varphi_n(y)\}_{n\in\mathbb N_0}$ satisfy the three--term recurrence relation
\begin{eqnarray}\label{ModMLPRec}
(n+2)\varphi_{n+1}(y) &=& 2y\varphi_{n}(y) - n\ \varphi_{n-1}(y) \quad (n\in \mathbb N)\\
\varphi_0(y)&=&2, \qquad \varphi_1(y)=2y.\nonumber
\end{eqnarray}
\end{thm}

\noindent{\it Proof.} The stated relation follows from the recurrence relation (\ref{MLPOrtho}).$\ \Box$

\begin{thm}
The generating function of sequence $\{\varphi_n(y)\}_{n\in\mathbb N_0}$ is given by
$$
\mathcal G(x,y)=\frac{\exp(2y\arctan x)-1}{xy}=\sum_{n=0}^\infty \varphi_n(y)x^n.
$$
\end{thm}

\noindent{\it Proof.} Starting from recurrence relation (\ref{ModMLPRec}) and summarizing, we have:
$$
\aligned
\sum_{n=1}^\infty (n+2)\varphi_{n+1}(y)x^n-2y\sum_{n=1}^\infty \varphi_n(y)x^n+\sum_{n=1}^\infty n\varphi_{n-1}(y)x^n&=0,\\
\sum_{n=2}^\infty (n+1)\varphi_{n}(y)x^{n-1}-2y\sum_{n=1}^\infty \varphi_n(y)x^n+\sum_{n=0}^\infty (n+1)\varphi_n(y)x^{n+1}&=0,\\
\frac1x\frac{\partial}{\partial x}\left(x\sum_{n=2}^\infty \varphi_{n}(y)x^{n}\right)-2y\sum_{n=1}^\infty \varphi_n(y)x^n
+x\frac{\partial}{\partial x}\left(x\sum_{n=0}^\infty \varphi_{n}(y)x^{n}\right)&=0,
\endaligned
$$
i.e.,
$$
(1+x^2)\frac{\partial}{\partial x}\big(x\mathcal G(x,y)\big)-2yx\mathcal G(x,y)-2=0.
$$
Solving obtained differential equation with initial condition \ $x\mathcal G(x,y)\big|_{x=0}=0$, we get generating function.$\ \Box$

Notice that
\begin{equation}\label{varphi(-y)}
\varphi_n(-y)=(-1)^n \ \varphi_n(y) \qquad (n\in\mathbb N).
\end{equation}

\begin{thm}
The polynomials of sequence $\{\varphi_n(y)\}_{n\in\mathbb N_0}$ satisfy
following orthogonality relation:
\begin{equation*}
\int_{-\infty}^{+\infty}\varphi_n(y)\varphi_m(y)\frac{y}{\sinh (\pi y)}\ dy = \frac{2}{n+1}\ \delta_{mn}\qquad(n,m\in\mathbb N_0).
\end{equation*}
\end{thm}

\noindent{\it Proof.} Using (\ref{MLPOrtho}) and (\ref{ModMLP}), we find
$$
-i^{n+m}\int_{-\infty}^{+\infty}\varphi_n(-y)\varphi_m(y)\frac{y^2}{y\sinh (\pi y)}\ dy = \frac{2}{n+1}\ \delta_{mn}\qquad(n,m\in\mathbb
N_0),
$$
what, with respect to (\ref{varphi(-y)}), gives orthogonality relation. $\ \Box$

\smallskip


The monic sequence
\begin{equation*}
\hat\varphi_n(y) = \frac{(n+1)!}{2^{n+1}}\ \varphi_n(y)  \qquad (n\in\mathbb N_0)
\end{equation*}
satisfies three term recurrence relation
\begin{eqnarray}\label{ModMLPRecMonic}
\hat\varphi_{n+1}(y) &=& y\hat\varphi_{n}(y) - \frac{n(n+1)}{4}\ \hat\varphi_{n-1}(y) \quad (n\in \mathbb N)\\
\hat\varphi_0(y)&=&1, \qquad \hat\varphi_1(x)=y.\nonumber
\end{eqnarray}

\begin{thm}
The exponential generating function of sequence $\{\hat\varphi_n(y)\}_{n\in\mathbb N_0}$ is given by
\begin{equation}\label{GenFun222}
\hat{\mathcal G}(x,y)=\frac{4\exp\bigl(2y\ \arctan (x/2)\bigr)}{x^2+4}=\sum_{n=0}^\infty \hat\varphi_{n}(y)\frac{x^n}{n!}.
\end{equation}
\end{thm}

\noindent{\it Proof.} In order to find exponential generating function of polynomials $\hat\varphi(y)$
\begin{equation*}\label{GenFun1}
\hat{\mathcal G}(x,y)=\sum_{n=0}^\infty \hat\varphi_{n}(y)\frac{x^n}{n!},
\end{equation*}
we will start with recurrence relation (\ref{ModMLPRecMonic}) and summarize
\begin{equation*}\label{GenFun2}
\sum_{n=1}^\infty \Bigl(\hat\varphi_{n+1}(y) - y\hat\varphi_{n}(y) + \frac{n(n+1)}{4}\ \hat\varphi_{n-1}(y)\Bigr)\frac{x^n}{n!} = 0,
\end{equation*}
i.e.
\begin{equation*}\label{GenFun3}
\sum_{n=2}^\infty \hat\varphi_{n}(y)\frac{x^{n-1}}{(n-1)!} - y\sum_{n=1}^\infty \hat\varphi_{n}(y)\frac{x^n}{n!} +
\frac{1}{4}\sum_{n=0}^\infty (n+1)(n+2)\ \hat\varphi_{n}(y)\frac{x^{n+1}}{(n+1)!} = 0.
\end{equation*}
Hence
\begin{equation*}\label{GenFun4}
\frac{\partial}{\partial x}\hat{\mathcal G}(x,y)-y\hat{\mathcal G}(x,y)+\frac14\frac{\partial}{\partial x}\Bigl(x^2\hat{\mathcal
G}(x,y)\Bigr)=0.
\end{equation*}
This is simple differential equation
$$
\frac{d\hat{\mathcal G}}{\hat{\mathcal G}}=2\frac{2y-x}{x^2+4}\ dx
$$
with initial value $\hat{\mathcal G}(0,y)=1$.
Its solution is the  function (\ref{GenFun222}). $\Box$

\begin{example}\rm
The first members of the sequence $\{\hat\varphi(y)\}_{n\in\mathbb N_0}$ are
$$
\hat\varphi_{0}(y) = 1,\quad
\hat\varphi_{1}(y) = y,\quad
\hat\varphi_{2}(y) = y^2-\frac12,\quad
\hat\varphi_{3}(y) = y^3-2y,
$$
$$
\hat\varphi_{4}(y) = y^4-5y^2+\frac32,\quad
\hat\varphi_{5}(y) = y^5-10y^3+\frac{23}2y.
$$
\end{example}

In the same manner, we can define {\it modified deformed Mittag-Leffler polynomials} by
\begin{equation*}
\varphi_n^{(h)}(y) = \frac{g_{n+1}^{(h)}(iy)}{i^{n+1}\ y}\qquad (n\in\mathbb N_0).
\end{equation*}
The properties of sequence $\{\varphi_n^{(h)}(y)\}_{n\in\mathbb N_0}$ can be derived in the same way as in case of
$\{\varphi_n(y)\}_{n\in\mathbb N_0}$. Hence, we give the results without proofs.

Firstly, the following is valid:
$$
\varphi_n^{(h)}(-y)=(-1)^n \ \varphi_n^{(h)}(y) \qquad (n\in\mathbb N).
$$

\begin{thm}
The successive members of sequence $\{\varphi_n^{(h)}(y)\}_{n\in\mathbb N_0}$ satisfy the three--term recurrence relation
\begin{eqnarray*}
(n+2)\varphi_{n+1}^{(h)}(y) &=& 2y\varphi_{n}^{(h)}(y) - h^2 n\ \varphi_{n-1}^{(h)}(y) \quad (n\in \mathbb N)\\
\varphi_0^{(h)}(y)&=&2, \qquad \varphi_1(y)^{(h)}=2y.\nonumber
\end{eqnarray*}
\end{thm}

\begin{thm}
The generating function of sequence $\{\varphi_n^{(h)}(y)\}_{n\in\mathbb N_0}$ is given by
$$
\mathcal G_h(x,y)=\frac{1}{xy}\left(\exp\left(2\frac{y}{h}\arctan hx\right)-1\right)=\sum_{n=0}^\infty \varphi_n^{(h)}(y)x^n.
$$
\end{thm}


\begin{thm}
The polynomials of sequence $\{\varphi_n^{(h)}(y)\}_{n\in\mathbb N_0}$ satisfy
following orthogonality relation:
\begin{equation*}
\int_{-\infty}^{+\infty}\varphi_n^{(h)}(y)\varphi_m^{(h)}(y)\frac{y}{\sinh (\pi y/h)}\ dy = \frac{2h^{2n}}{n+1}\
\delta_{mn}\qquad(n,m\in\mathbb N_0).
\end{equation*}
\end{thm}

The monic sequence
\begin{equation*}
\hat\varphi_n^{(h)}(y) = \frac{(n+1)!}{2^{n+1}}\ \varphi_n^{(h)}(y)  \qquad (n\in\mathbb N_0)
\end{equation*}
satisfies three term recurrence relation
\begin{eqnarray*}
\hat\varphi_{n+1}^{(h)}(y) &=& y\hat\varphi_{n}^{(h)}(y) - \frac{h^2}{4}n(n+1)\ \hat\varphi_{n-1}^{(h)}(y) \quad (n\in \mathbb N)\\
\hat\varphi_0^{(h)}(y)&=&1, \qquad \hat\varphi_1^{(h)}(x)=y.\nonumber
\end{eqnarray*}

\begin{thm}
The exponential generating function of $\{\hat\varphi_n^{(h)}(y)\}_{n\in\mathbb N_0}$ is
$$
\hat{\mathcal G}_h(x,y)=(4+h^2 x^2)^{-1/h^2}\exp\Bigl(2\frac yh \arctan \frac {hx}{2}\Bigr)=\sum_{n=0}^\infty
\hat\varphi_{n}^{(h)}(y)\frac{x^n}{n!}\ .
$$
\end{thm}

\begin{example}\rm
The first members of the sequence $\{\hat\varphi^{(h)}(y)\}_{n\in\mathbb N_0}$ are
$$
\hat\varphi_{0}^{(h)}(y) = 1,\quad
\hat\varphi_{1}^{(h)}(y) = y,\quad
\hat\varphi_{2}^{(h)}(y) = y^2-\frac12,\quad
\hat\varphi_{3}^{(h)}(y) = y^3-2y,
$$
$$
\hat\varphi_{4}^{(h)}(y) = y^4-5y^2+\frac32,\quad
\hat\varphi_{5}^{(h)}(y) = y^5-10y^3+\frac{23}2y.
$$
\end{example}

The modified Mittag-Leffler $\{\varphi_n(y)\}$ and deformed Mittag--Leffler polynomials $\{\varphi_n^{(h)}(y)\}$
are real polynomials and, because of orthogonality, they have all real zeros.

\end{document}